\documentclass[12pt]{amsart}

\usepackage{amsmath}
\usepackage[curve]{xypic}
\usepackage{caption}
\usepackage{xspace}
\usepackage{cite}
\usepackage{enumitem}
\usepackage{color}
\usepackage{url}
\usepackage[usenames,dvipsnames]{xcolor}
\usepackage{graphicx}

\usepackage{hyperref}%[hidelinks]

\makeatletter % makes "@" a normal character so that it can be used in macros
\def\@cite#1#2{{\m@th\upshape\bfseries%
[{#1\if@tempswa{\m@th\upshape\mdseries, #2}\fi}]}}
\makeatother % returns "@" to its normal state

\theoremstyle{plain}
\newtheorem{thm}{Theorem}[section]% subsection

\newtheorem{lem}[thm]{Lemma}

\theoremstyle{definition}

\newtheorem{ex}[thm]{Example}

\theoremstyle{remark}

\numberwithin{equation}{subsection}
\renewcommand{\bold}[1]{\medskip \noindent {\bf #1 }\nopagebreak}

\newcommand{\nc}{\newcommand}
\newcommand{\rnc}{\renewcommand}
\nc\bA{\mathbb{A}}
\nc\bB{\mathbb{B}}
\nc\bC{\mathbb{C}}
\nc\bD{\mathbb{D}}
\nc\bE{\mathbb{E}}
\nc\bF{\mathbb{F}}
\nc\bG{\mathbb{G}}
\nc\bH{\mathbb{H}}
\nc\bI{\mathbb{I}}
\nc{\bJ}{\mathbb{J}} 
\nc\bK{\mathbb{K}}
\nc\bL{\mathbb{L}}
\nc\bM{\mathbb{M}}
\nc\bN{\mathbb{N}}
\nc\bO{\mathbb{O}}
\nc\bP{\mathbb{P}}
\nc\bQ{\mathbb{Q}}
\nc\bR{\mathbb{R}}
\nc\bS{\mathbb{S}}
\nc\bT{\mathbb{T}}
\nc\bU{\mathbb{U}}
\nc\bV{\mathbb{V}}
\nc\bW{\mathbb{W}}
\nc\bY{\mathbb{Y}}
\nc\bX{\mathbb{X}}
\nc\bZ{\mathbb{Z}}
\nc\cA{\mathcal{A}}
\nc\cB{\mathcal{B}}
\nc\cC{\mathcal{C}}
\rnc\cD{\mathcal{D}}
\nc\cE{\mathcal{E}}
\nc\cF{\mathcal{F}}
\nc\cG{\mathcal{G}}
\rnc\cH{\mathcal{H}}
\nc\cI{\mathcal{I}}
\nc{\cJ}{\mathcal{J}} 
\nc\cK{\mathcal{K}}
\rnc\cL{\mathcal{L}}
\nc\cM{\mathcal{M}}
\nc\cN{\mathcal{N}}
\nc\cO{\mathcal{O}}
\nc\cP{\mathcal{P}}
\nc\cQ{\mathcal{Q}}
\rnc\cR{\mathcal{R}}
\nc\cS{\mathcal{S}}
\nc\cT{\mathcal{T}}
\nc\cU{\mathcal{U}}
\nc\cV{\mathcal{V}}
\nc\cW{\mathcal{W}}
\nc\cY{\mathcal{Y}}
\nc\cX{\mathcal{X}}
\nc\cZ{\mathcal{Z}}

\nc{\dmo}{\DeclareMathOperator}
\dmo{\Tw}{Twist}
\dmo{\CP}{Pres}
\rnc{\Re}{\operatorname{Re}}
\rnc{\Im}{\operatorname{Im}}
\rnc{\span}{\operatorname{span}}
\dmo{\rank}{rank}
\dmo{\End}{End}
\dmo{\Hom}{Hom}
\dmo{\Jac}{Jac}
\dmo{\Id}{Id}
\dmo{\lcm}{lcm}
\dmo{\Area}{Area}

\nc{\Tm}{Teichm\"uller\xspace}
\nc{\odd}{\cH^{\text{odd}}(4)}
\nc{\hyp}{\cH^{\text{hyp}}(4)}
\nc{\prym}{\tilde{\mathcal{Q}}(3,-1^3)}
\nc{\G}{GL^+(2,\bR)}
\nc{\GL}{GL^+}

\begin{document}
%%%%%%%%%%%%%%%%%%%%%%%%%%%%%%%%%%%%%%%%%%%

\title{From rational billiards to dynamics on moduli spaces}
\author[Wright]{Alex~Wright}
%\address{\hspace{-0.5cm} Math\ Department\newline
%University of Chicago\newline
%5734 South University Avenue\newline
%Chicago, IL 60637}
\email{alexmurraywright@gmail.com}
%
%\subjclass[2010]{22E60, 15A57, 17B20, 58C35}
%\keywords{}
%\date{\today}
%\dedicatory{Preliminary version. Comments welcome.}

\begin{abstract}
%Subtitle: The work of Eskin and Mirzakhani draws inspiration from homogeneous spaces to study more complicated spaces that arise in diverse areas of mathematics and physics, with applications to some deceptively simple dynamical systems.
This short expository note  gives an elementary introduction to the study of dynamics on certain moduli spaces, and in particular the recent breakthrough result of Eskin, Mirzakhani, and Mohammadi. 
We also discuss the context and applications of this result, and connections to other areas of mathematics such as algebraic geometry, Teichm\"uller theory, and ergodic theory on homogeneous spaces.

%%Starting with the problem of billiards in a polygon, we explain that under a rationality assumption the polygon can be ``unfolded" to give a surface with extra structure, and there is an $SL(2,\mathbb{R})$ action on the space of all such surfaces. 

%Consider a billiard ball bouncing around in a polygon. This simple system demonstrates remarkable complexity. For example, it is an open problem to prove that there is a periodic billiard trajectory in every polygon. However, if the angles are all rational multiples of $\pi$, a great deal is known. This is because such a polygon can be ``unfolded" to give a surface with extra structure, and there is an $SL(2,\mathbb{R})$ action on the space of all such surfaces. We will explain the relevance of this action, and state a recent result of Eskin, Mirzakhani, and Mohammadi, which gives that the closure of every $SL(2,\mathbb{R})$ orbit is a manifold. Applications and connections to other areas of mathematics will be discussed. 
\end{abstract}

%\date{\today.}
%\dedicatory{Preliminary version. Comments welcome.}

\maketitle
% removes page number from first page
\thispagestyle{empty}

%\newpage

%Five Fields Medalists have studied the GL(2,R) action on the moduli space of translation surfaces. Come learn the background, and the newest breakthrough. 
%
%
%The study of the GL(2,R) action on the moduli space of translation surfaces, sometimes called Teichmuller dynamics, is a thriving field rich with connections to diverse areas of mathematics. Work in this area was mentioned in the citations for two of the four fields medals awarded in 2014: those of Artur Avila and Maryam Mirzakhani. Here we will focus on the work of Eskin, Mirzakhani, and Mohammadi. For a much broader survey of MirzakhaniÕs work, see the survey of McMullen. 

%%%%%%%%%%%%%%%%%%% 
% TABLE OF CONTENTS
%%%%%%%%%%%%%%%%%%%
% allows subsections (depth 1) to be displayed in table of contents
\setcounter{tocdepth}{1} 
\tableofcontents
%\vfill
%\newpage

\section{Rational billiards} 
Consider a point bouncing around in a polygon. Away from the edges, the point moves at unit speed. At the edges, the point bounces according to the usual rule that angle of incidence equals angle of reflection. If the point hits a vertex, it stops moving. The path of the point is called a billiard trajectory.   

The study of billiard trajectories is a basic problem in dynamical systems and arises naturally in physics. For example, consider two points of different masses moving on a interval, making elastic collisions with each other and with the endpoints. This system can be modeled by billiard trajectories in a right angled triangle \cite{MT}. 

A rational polygon is a polygon all of whose angles are rational multiples of $\pi$. Many mathematicians are especially interested in billiards in rational polygons for the following three reasons.

First, without the rationality assumption,  few tools are available, and not much is known. For example, it is not even known if every triangle has a periodic billiard trajectory. With the rationality assumption, quite a lot can be proven.

Second, even with the rationality assumption a wide range of interesting behavior is possible, depending on the choice of polygon.  

Third,  the rationality assumption leads to surprising and beautiful connections to algebraic geometry, Teichm\"uller theory, ergodic theory on homogenous spaces, and other areas of mathematics. 
 
The assumption of rationality first arose from the following simple thought experiment. What if, instead of letting a billiard trajectory bounce off an edge of a polygon, we instead allowed the trajectory to continue straight, into a reflected copy of the polygon? 

\begin{figure}[h]
\includegraphics[width=0.35\linewidth]{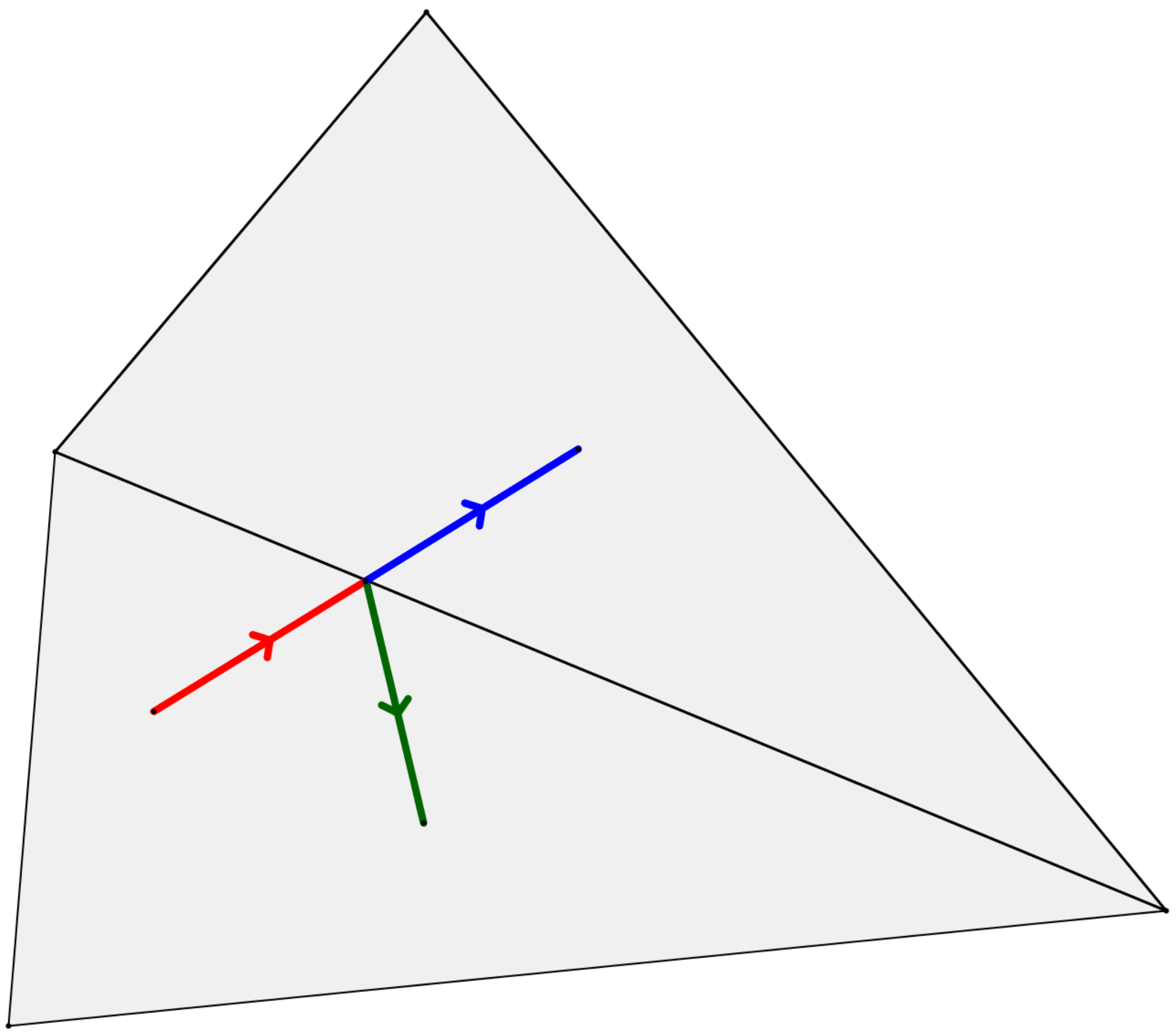}
\caption{A billiard trajectory in a polygon $P$. Instead of allowing the trajectory to bounce off the edge of $P$, we may allow it to continue straight into a reflected copy of $P$. A key observation is that the trajectory that continues into the reflected copy of $P$ is in fact the reflection of the trajectory in $P$ that bounces off of the edge. }
\label{F:Reflect}
\end{figure}

This leads us to define the ``unfolding" of a polygon $P$ as follows: Let $G$ be the subgroup of $O(2)$ (linear isometries of $\bR^2$) generated by the derivatives of reflections in the sides of $P$. The group $G$ is finite if and only if the polygon $P$ is rational (in which case $G$ is a dihedral group). For each $g\in G$, consider the polygon $gP$. These polygons $gP$ can be translated so that they are all disjoint in the plane. We identify the edges in pairs in the following way. Suppose $r$ is the derivative of the reflection in one of the edges of $hP$. Then this edge of $hP$ is identified with the corresponding edge of $rhP$. 

\begin{figure}[ht]
\begin{minipage}[c]{0.63\linewidth}
\centering
\includegraphics[width=\textwidth]{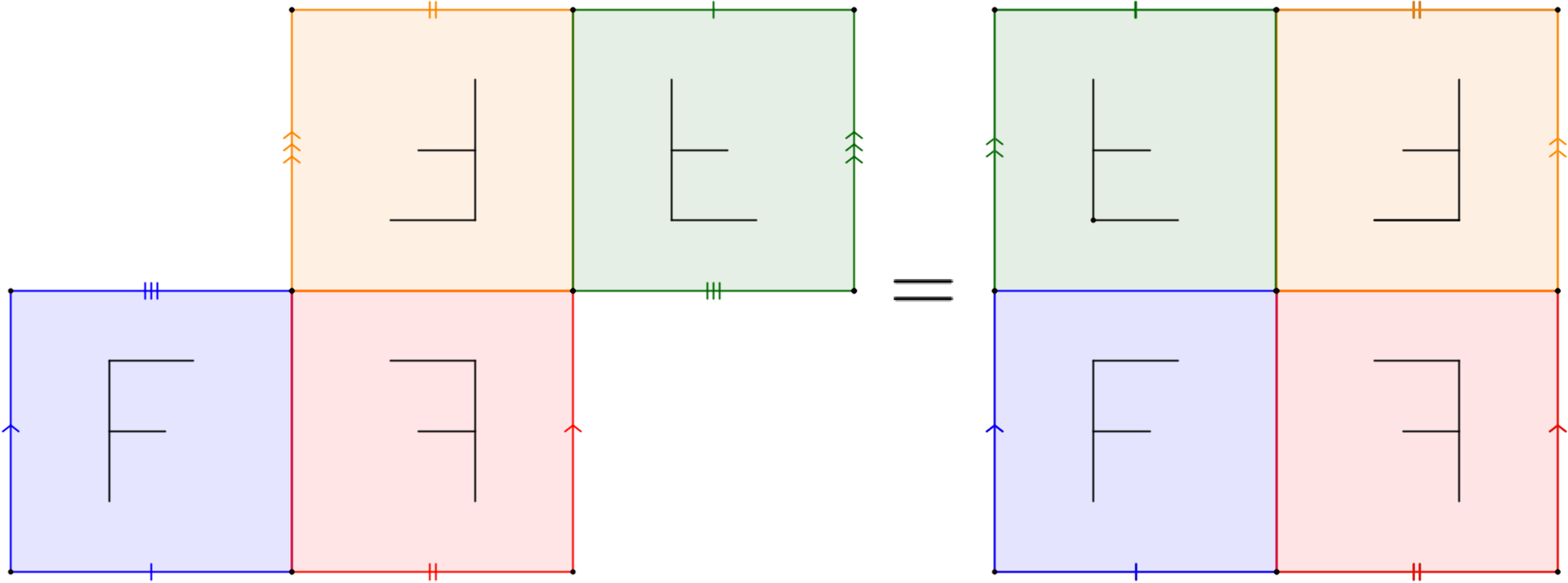}
\end{minipage}
\hspace{0.5cm}
\begin{minipage}[c]{0.27\linewidth}
\centering
\includegraphics[width=\textwidth]{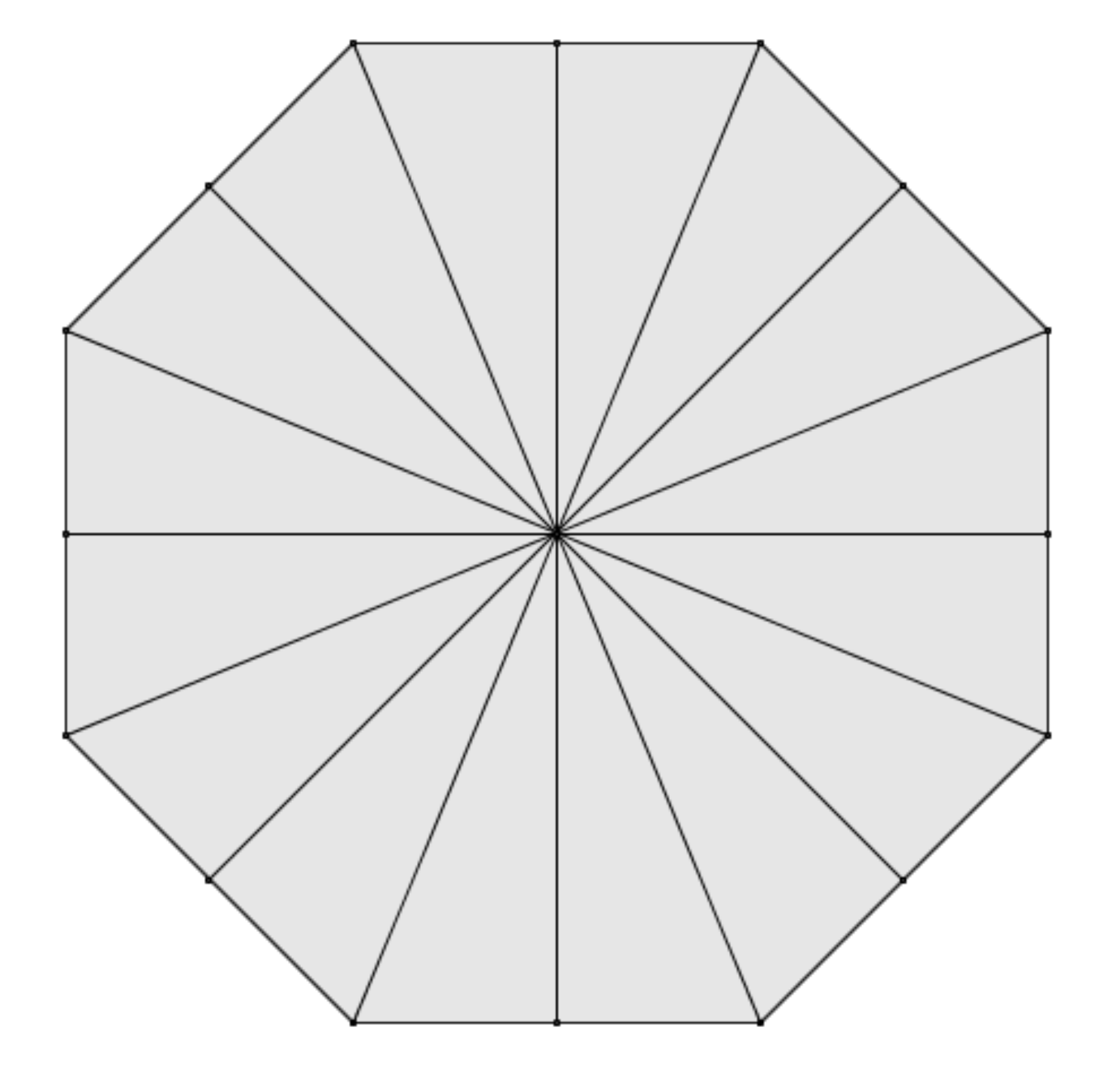}
\end{minipage}
\caption{Left: The unit square unfolds to four squares, with opposite edges identified (a flat torus). (By ``opposite edges" in these pictures, we mean pairs of boundary edges that are perpendicularly across from each other, so for example the top left and top right vertical edges on the unfolding of the square are opposite.) When two polygons are drawn with an adjacent edge, by convention this means these two adjacent edges are identified. Here each square has been decorated with the letter F, to illustrate which squares are reflections of other squares.   Right: Unfolding the right angled triangle with smallest angle $\pi/8$ gives the regular octagon with opposite sides identified.}
\label{F:Unfold}
\end{figure}
\begin{figure}[h!]
\centering
\includegraphics[scale=0.2]{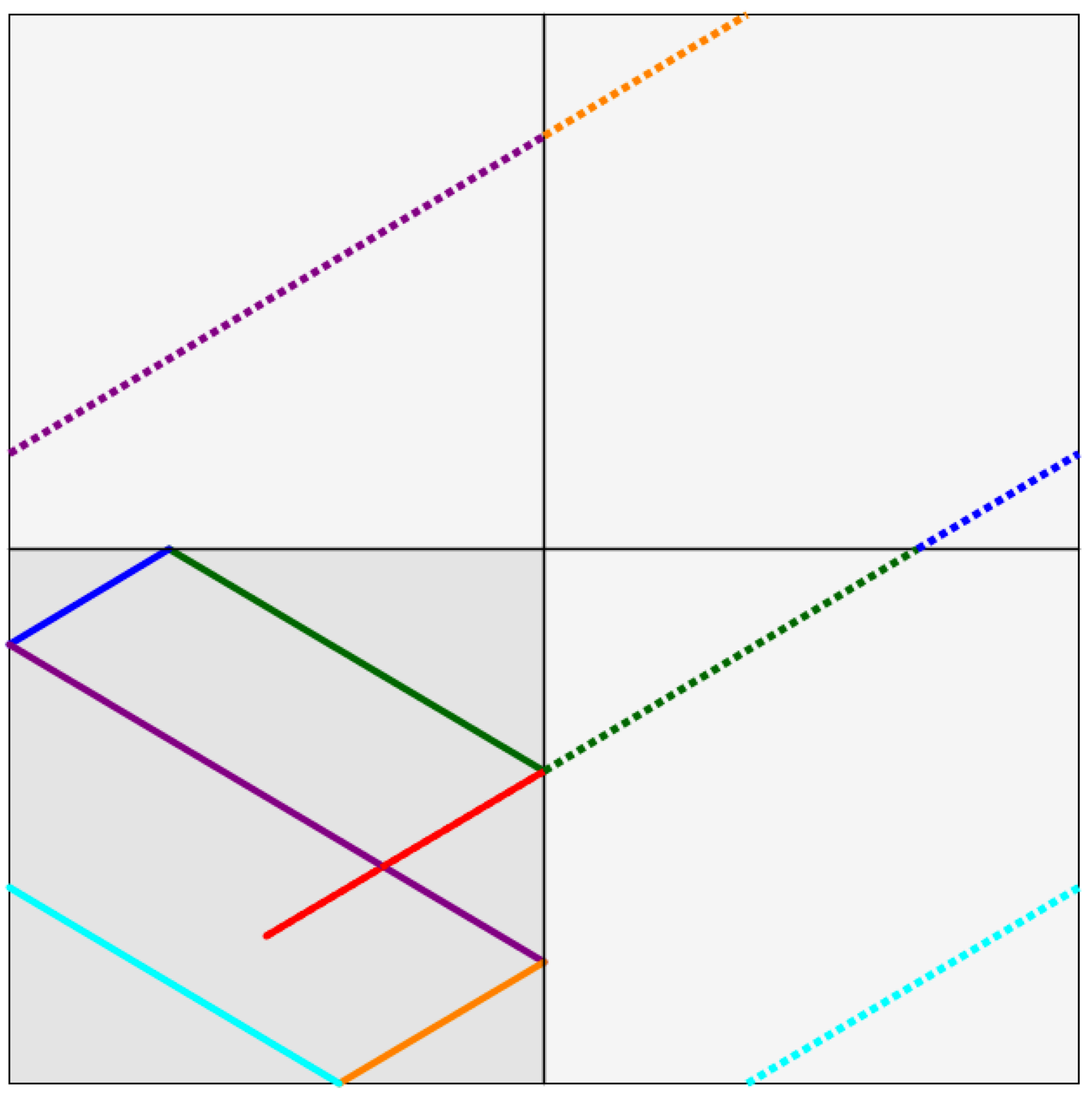}
\caption{A billiard trajectory on a rational polygon unfolds to a straight line on the unfolding of the polygon. In this illustration, we have unfolded a billiard trajectory on square (bottom left) to a straight line on a flat torus. The square and its unfolding are superimposed, the billiard trajectory is drawn with a solid line, and the unfolded straight line is drawn with a dotted line.}
\label{F:CurrentEvents}
\end{figure}
The unfolding construction is most easily understood through examples: see Figures \ref{F:Unfold} and \ref{F:CurrentEvents}.

\section{Translation surfaces}

Unfoldings of rational polygons are special examples of translation surfaces. There are several equivalent definitions of translation surface, the most elementary of which is a finite union of polygons in in the plane with edge identifications, obeying certain rules, up to a certain equivalence relation. The rules are: 
\begin{enumerate}
\item The interiors of the polygons must be disjoint, and if two edges overlap then they must  be identified. 
\item Each edge is identified with exactly one other edge, which must be a translation of the first. The identification is via this translation. 
\item When an edge of one polygon is identified with an edge of a different polygon, the polygons must be on ``different sides" of the edge. For example, if a pair of vertical edges are identified, one must be on the left of one of the polygons, and the other must be on the right of the other polygon. 
\end{enumerate}
 Two such families of polygons are considered to be equivalent if they can be related via a string of the following ``cut and paste" moves. 
\begin{enumerate}
\item A polygon can be translated. 
\item A polygon can be cut in two along a straight line, to give two adjacent polygons. 
\item Two adjacent polygons that share an edge can be glued to form a single polygon. 
\end{enumerate}

\begin{figure}[h!]
\centering
\includegraphics[width=\textwidth]{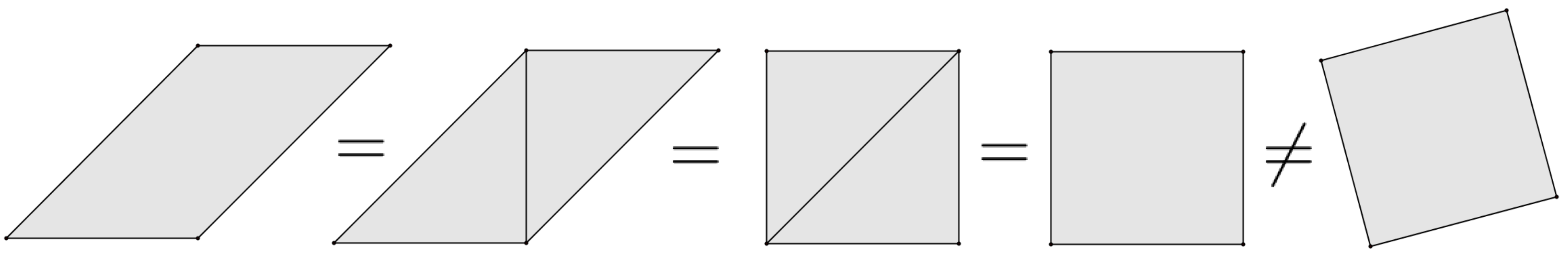}
\caption{In all five translation surfaces above, opposite edges are identified. In the leftmost four, each adjacent pair of translation surfaces differs by one of the above three moves, so all four of these pictures give the same translation surface. The rightmost rotated surface is (presumably) not equal to the other four, since rotation is not one of the three allowed moves.}
\label{F:CutAndPaste}
\end{figure}

In general, it is difficult to decide if  two collections of polygons as above are equivalent (describe the same translation surface), because each collection of polygons is equivalent to infinitely many others. 

\begin{figure}[h!]
\centering
\includegraphics[width=0.5\textwidth]{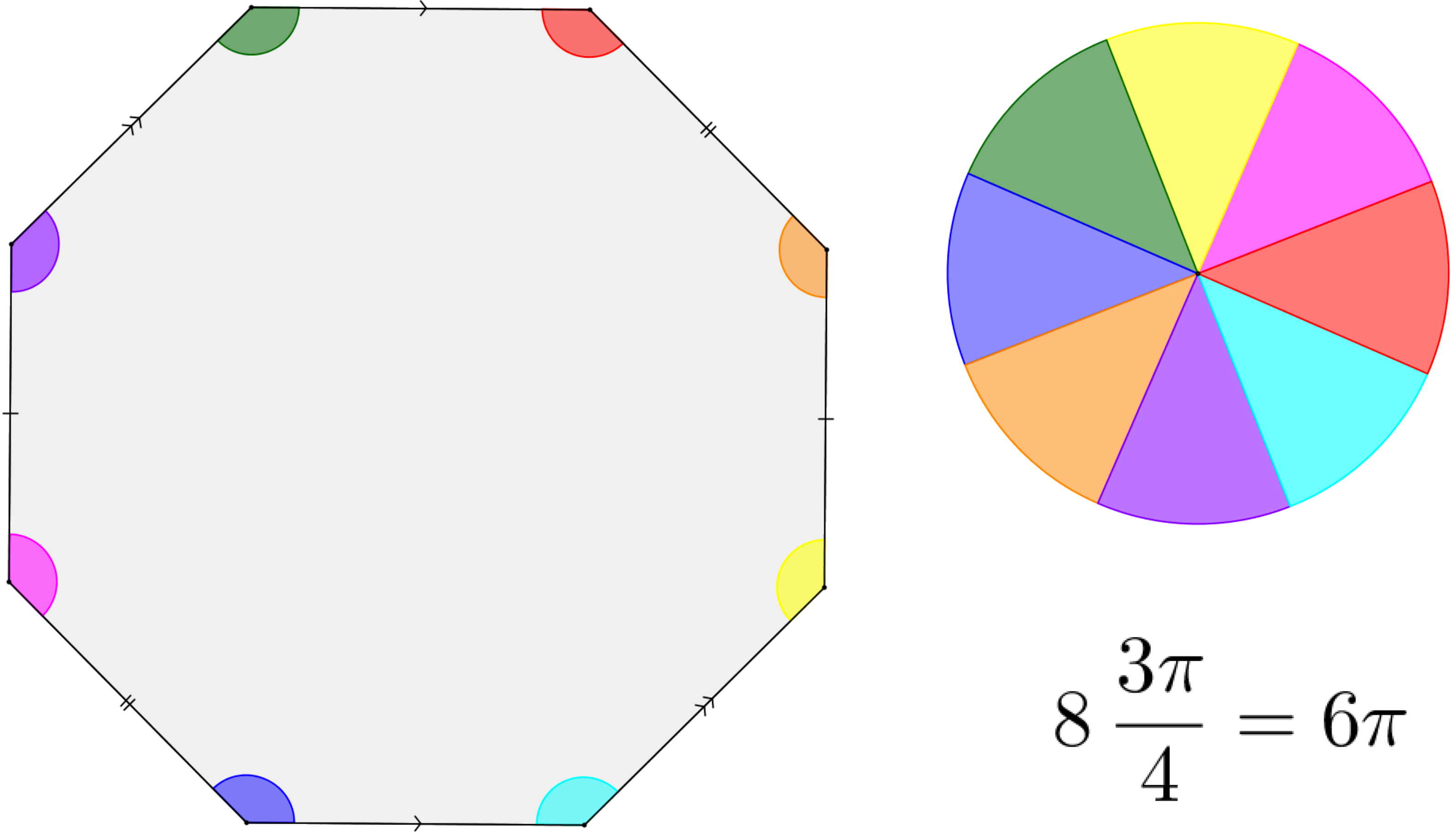}
\caption{The edge identifications imply that the 8 corners of the octagon are in fact all identified, and give a single point on the translation surface. Around this point there is $6\pi$ total angle, since at each of the 8 corners of the polygon there is $\frac34\pi$ interior angle.}
\label{F:ConeAngle}
\end{figure}

The requirements above ensure that the union of the polygons modulo edge identifications gives a closed surface. This surface has flat metric, given by the flat metric on the plane, away from a finite number of singularities. The singularities arise from the corners of the polygons. For example, in the regular 8-gon with opposite sides identified, the 8 vertices give rise to a single point with cone angle. See Figure \ref{F:ConeAngle}.

The singularities of the flat metric on a translation surface are always of a very similar conical form, and the total angle around a singularity on a translation surface is always an integral multiple of $2\pi$. Note that, although the flat metric is singular at these points, the underlying topological surface is not singular at any point. (That is, at every single point, including the singularities of the flat metric, the surface is locally homeomorphic to $\bR^2$.) 

Most translation surfaces do not arise from unfoldings of rational polygons. This is because unfoldings of polygons are exceptionally symmetric, in that they are tiled by isometric copies of the polygon.  

\begin{figure}[h!]
\centering
\includegraphics[scale=0.4]{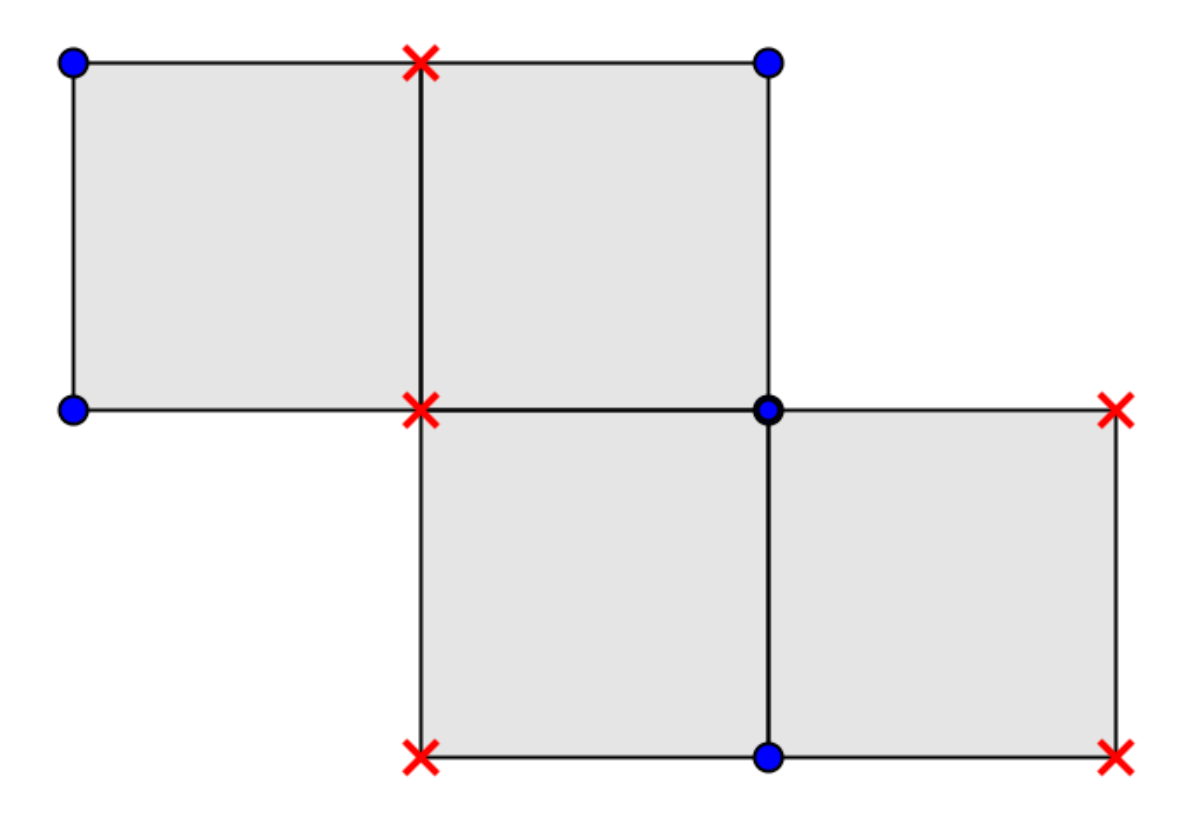}
\caption{Consider the translation surface described by the above polygon, with opposite edges identified. This surface has two singularities,  each with total angle $4\pi$. One singularity has been labelled with a dot, and the other with an x. An Euler characteristic calculation ($V-E+F=2-2g$) show that it has genus 2. The regular octagon with opposite sides identified also has genus 2, but it has only a single singularity, with total angle $6\pi$.}
\label{F:STgenus2}
\end{figure}

Translation surfaces satisfy a Gauss-Bonnet type theorem. If a translation surface has $s$ singularities with cone angles $$(1+k_1)2\pi, (1+k_2)2\pi, \ldots, (1+k_s)2\pi,$$ then the genus $g$ is given by the formula $2g-2=\sum k_i$. (So, in a formal comparison to the usual Gauss-Bonnet formula, one might say that each extra $2\pi$ of angle on a translation surface counts for 1 unit of negative curvature.)

Consider now the question of how a given translation surface can be deformed to give other translation surfaces. The polygons, up to translation, can be recorded by their edge vectors in $\bC$ (plus some finite amount of combinatorial data, for example the cyclic order of edges around the polygons). Not all edge vectors need be recorded, since some are determined by the rest. Changing the edge vectors (subject to the condition that identified edges should remain parallel and of the same length) gives a deformation of the translation surface. 
\begin{figure}[h!]
\centering
\includegraphics[width=\linewidth]{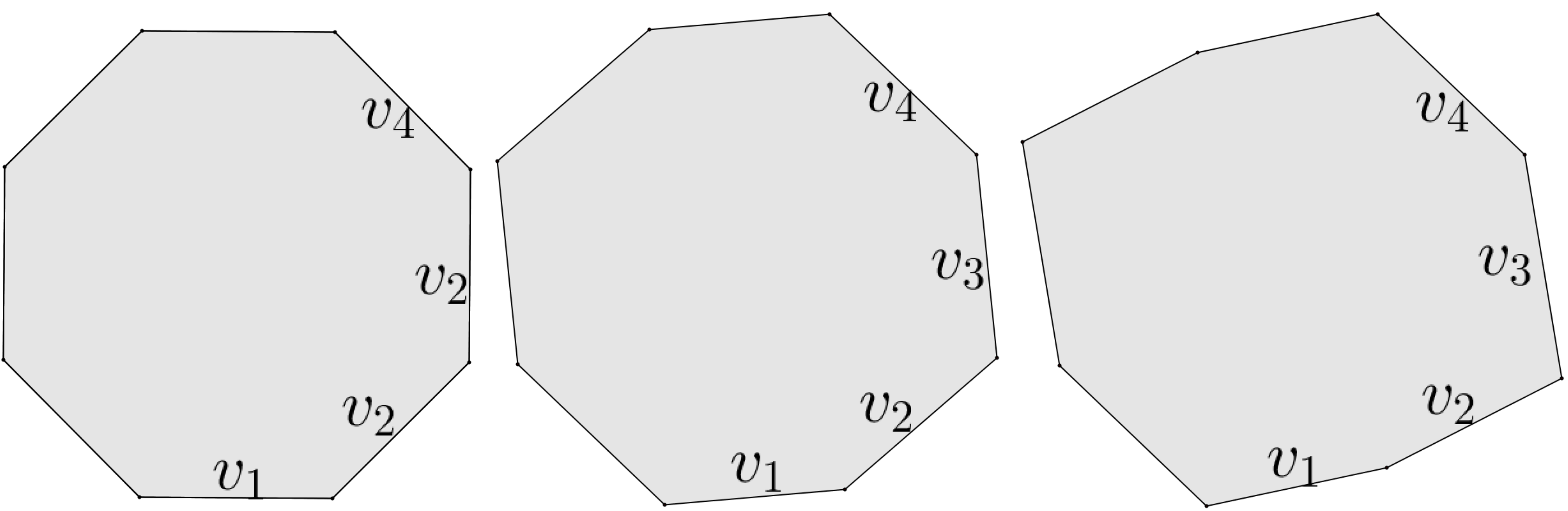}
\caption{Any octagon whose opposite edges are parallel can be described by the 4-tuple of its edges vectors $(v_1,v_2,v_3,v_4)\in \bC^4$. (Not all choices of $v_i$ give valid octagons.) The coordinates $(v_1, v_2, v_3, v_4)$ are local coordinates for space of deformations of the regular octagon translation surface. These coordinates are not canonical: other equally good coordinates can be obtained by cutting up the octagon and keeping track of different edge vectors. }
\label{F:STgenus2}
\end{figure}

To formalize this observation, we define moduli spaces of translation surfaces. Given an unordered collection $k_1, \ldots, k_s$ of positive integers whose sum is $2g-2$, the stratum $\cH(k_1, \ldots, k_s)$ is defined to be the set of all translation surfaces with $s$ singularities, of cone angles  $(1+k_i)2\pi, i=1,\ldots, s$. The genus of these surfaces must be $g$ by the Gauss-Bonnet formula above. We have

\begin{lem}
Each stratum is a complex orbifold of dimension $n=2g+s-1$. Each stratum has a finite cover that is a manifold and has an atlas of charts to $\bC^n$ with transition functions in $GL(n,\bZ)$. 
\end{lem}

\begin{figure}[h!]
\centering
\includegraphics[width=0.7\linewidth]{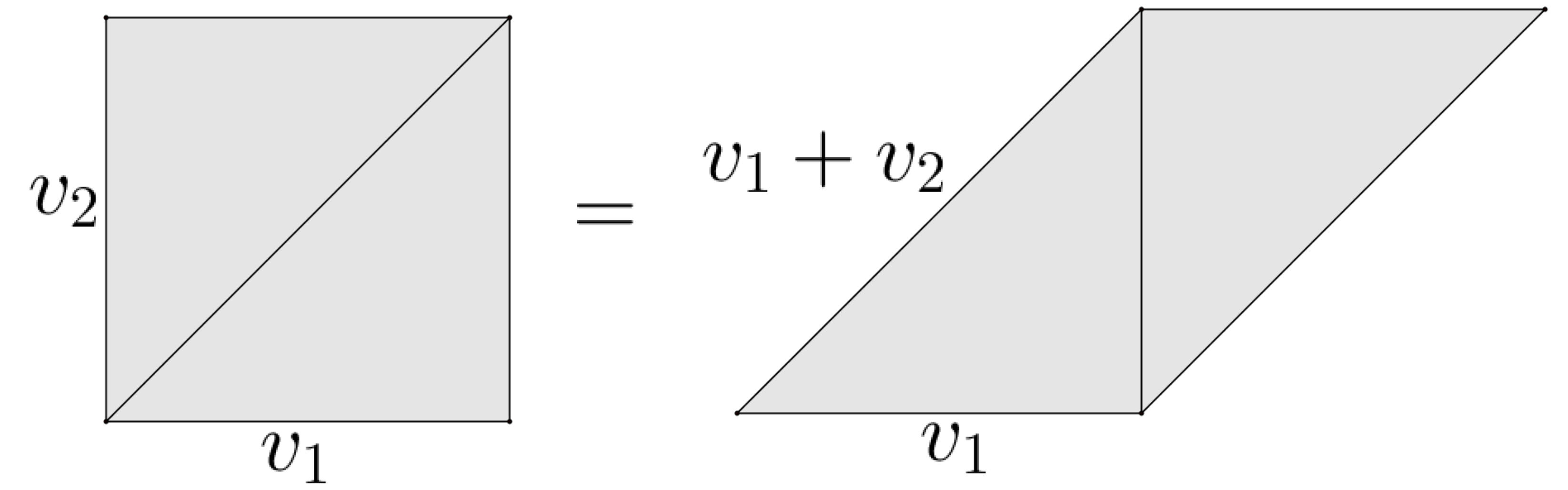}
\caption{These two polygons (with opposite sides identified) both describe the same translation surface. Keeping track of the edge vectors in either polygon gives equally good local coordinates for the space of nearby translation surfaces. The two local coordinates thus obtained are related by the linear transformation $(v_1, v_2)\mapsto (v_1, v_1+v_2)$. }
\label{F:NewCoords}
\end{figure}

The coordinate charts are called period coordinates. They consist of complex edge vectors of polygons.  That strata are orbifolds instead of manifolds is a technical point that should be ignored by non-experts. 

 Strata are not always connected, but their connected components have been classified by Kontsevich and Zorich \cite{KZ}. There are always at most three connected components. The topology (and  birational geometry) of strata is currently not well understood. Kontsevich has conjectured that strata are $K(\pi,1)$ spaces. 

\section{The $GL(2, \bR)$ action}

There is a $GL(2, \bR)$ action on each stratum, obtained by acting linearly on polygons and keeping the same identification. 

\begin{figure}[h!]
\centering
\includegraphics[width=0.7\textwidth]{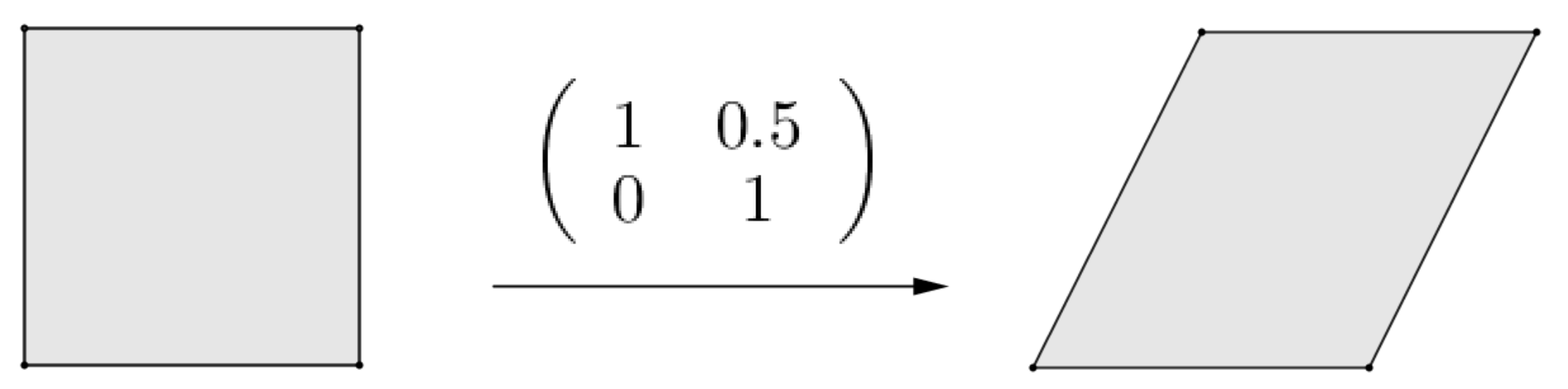}
\caption{An example of the $GL(2, \bR)$ action. In both pictures, opposite edges are identified. }
\label{F:Shear}
\end{figure}
Note that if two edges or polygons differ by translation by a vector $v$, then their images under the linear map $g\in GL(2,\bR)$ must  differ by translation by $gv$.

\begin{ex}
The stabilizer of the standard flat torus (a unit square with opposite sides identified) is $GL(2,\bZ)$. For example, Figure \ref{F:CutAndPaste} (near the beginning of the previous section) proves that $\left(\begin{array}{cc}1&1\\0&1\end{array}\right)$ is in the stabilizer. This example illustrates the complexity of the $GL(2,\bR)$ action: applying a large matrix (say of determinant 1) will yield a collection of very long and thin polygons, but it is hard to know when this collection of  polygons is equivalent to a more reasonable one. 
\end{ex}

Translation surfaces have a well defined area, given by the sum of the areas of the polygons. The action of $SL(2,\bR)$ of determinant 1 matrices in $GL(2,\bR)$ preserves the locus of unit area translation surfaces. This locus is not compact,  because the polygons can have edges of length going to 0, even while the total area stays constant.

Define $$g_t=\left(\begin{array}{cc} e^t&0\\0&e^{-t}\end{array}\right)\in SL(2,\bR).$$ Suppose one wants to know if a translation surface $S$ has a vertical line joining singularities (cone points) of length $e^{10}$. This is equivalent to the question of whether $g_{10}(S)$ has a vertical line segment of length 1 joining two singularities. 

\begin{figure}[h!]
\centering
\includegraphics[width=0.7\textwidth]{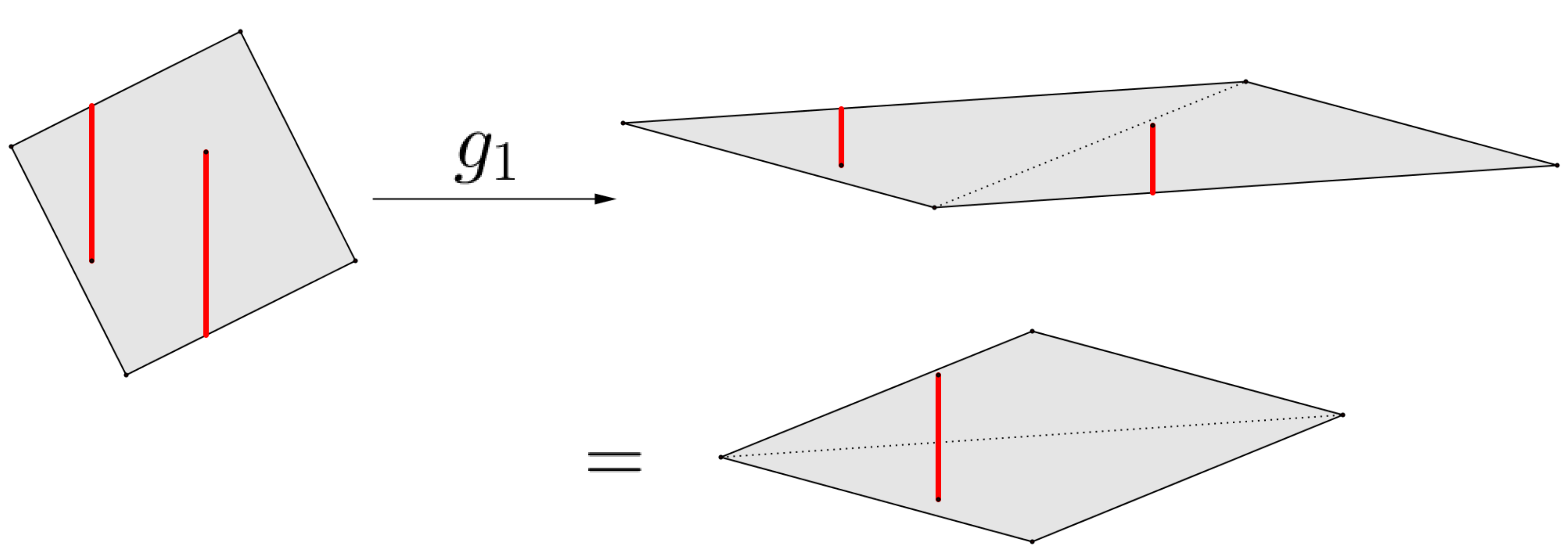}
\caption{Applying $g_t$ to a translation surface with a vertical line segment gives gives a new translation surface with a shorter vertical line segment.}
\label{F:Shear}
\end{figure}

In fact, for any matrix $g\in GL(2,\bR)$, the surfaces $S$ and $g(S)$ are closely related, since one can go back and forth between them using $g$ and $g^{-1}$. Really $S$ and $g(S)$  are just different perspectives on the same object, in which different features are apparent. To understand $S$ from all possible perspectives, we'd like to understand its $GL(2,\bR)$ orbit. However, just from definitions, its not really possible to understand the $GL(2,\bR)$ orbit of a surface. It's really hard to know, given two surfaces $S$ and $S'$, if there are large matrices $g$ so that the polygons defining $g(S)$ can be cut up and reglued to be almost equal to the polygons defining $S'$!

\section{Renormalization}

This section aims to give some of the early motivations and successes of the field. It can be safely skipped and returned to later by anyone eager to get to the modern breakthrough and its applications and connections to other area of mathematics. 

 Suppose once again that we have a vertical line segment of length $L$ on a translation surface $S$. For example, if $S$ is the unfolding of a rational polygon, the vertical line might be the unfolding of a billiard trajectory. If one is interested in a line that is not vertical, one can rotate the whole picture (giving a different translation surface) so that it becomes vertical. 

Applying $g_t$  results  in a translation surface $g_t(S)$ with a vertical segment of length $e^{-t}L$. We are interested in doing this when $L$ is very large, and $t=\log(L)$ is chosen so the new vertical segment will have length 1. Indeed the point is to do this over and over as $L$ gets longer, giving a family of surfaces $g_t(S)$.

This idea of taking longer and longer trajectories (here vertical lines on the translation surface) and replacing them with bounded length trajectories on  new objects is called renormalization, and is a powerful and frequently used tool in the study of dynamical systems. The typical strategy is to transfer some understanding of the sequence of renormalized objects into results on the behavior of the original system. In this case, showing that the geometry of $g_t(S)$ does not degenerate allows good understanding of vertical lines on $S$. 

\begin{thm}[Masur's criterion \cite{Ma}]
Suppose $\{g_t(S): t\geq 0\}$ does not diverge to infinity in the stratum. Then every infinite vertical line on $S$ is equidistributed on $S$. 
\end{thm}

``Equidistributed" is a technical term that indicates that the vertical lines becomes dense in $S$ without favoring one part of $S$ over another. Using this, Kerckhoff-Masur-Smillie \cite{KMS} were able to show

\begin{thm}
In every translation surface, for almost every slope, every infinite line of this slope is equidistributed. 
\end{thm} 

There are some surfaces where much more is true. For example, on the unit square with opposite sides identified, any line of rational slope is periodic, and every line of irrational slope is equidistributed. Genus one translation surfaces are quite special, because $GL(2,\bR)$ acts transitively on the space of genus one translation surfaces. In particular, the $GL(2,\bR)$ orbit of any genus one translation surface is closed, in a trivial way, since the orbit is the entire moduli space. 

\begin{thm}[Veech Dichotomy]
If $S$ is a translation surface with closed $GL(2,\bR)$ orbit, then for all but countably many slopes, every line with that slope is equidistributed. Moreover every line with slope contained in the countable set is periodic. 
\end{thm}

\begin{figure}[h!]
\centering
\includegraphics[scale=0.3]{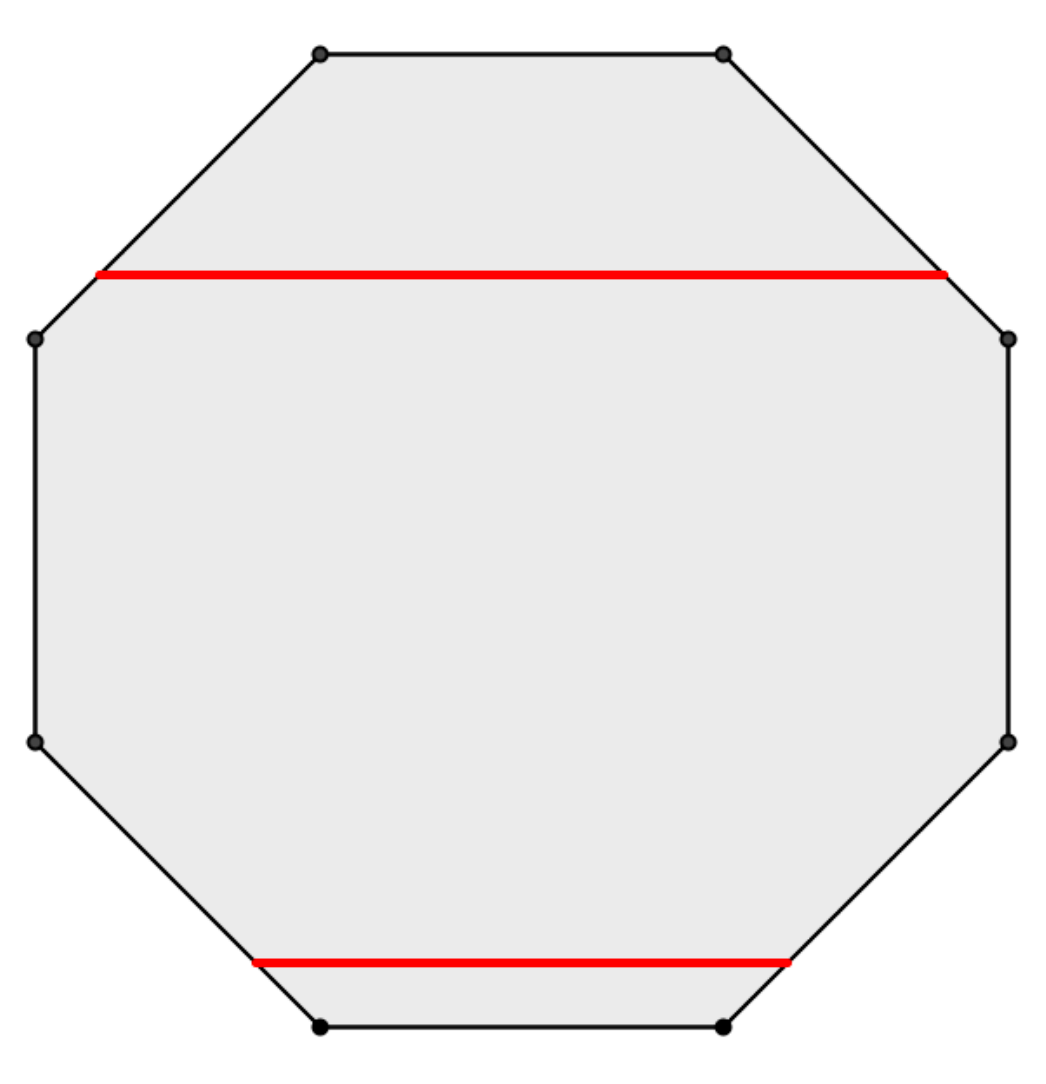}
\caption{An example of a periodic line on the regular octagon with opposite sides identified.}
\label{F:Periodic}
\end{figure}

Veech also showed that the regular $2n$-gon with opposite sides identified has closed orbit. However, the property of having a closed orbit is extremely special. 

\begin{thm}[Masur \cite{Ma2}, Veech \cite{V2}]
The $GL(2,\bR)$ orbit of almost every translation surface is dense in a connected component of a stratum. 
\end{thm}

For the experts, we remark that in fact Masur and Veech showed the stronger statement that the $g_t$ action on the loci of unit area surfaces in a connected component of a stratum is ergodic, with respect to a Lebesgue class probability measure called the Masur-Veech measure. 

This result of Masur and Veech is not satisfactory from the point of view of billiards in rational polygons, since the set of translation surfaces that are unfoldings of polygons has measure 0.  

\section{Eskin-Mirzakhani-Mohammadi's breakthrough}

%In light of the idea of renormalization, it is important to understand the $SL(2,\bR)$ orbit closures of all translation surfaces. It is equivalent to understand the $GL(2,\bR)$ orbit closures, since the difference is just scaling of the area. 
%
%In fact it would also be helpful to know more specific information, such as the $g_t$ orbit closures. However,  

The $g_t$ orbit closure of a translation surface may be a fractal object. While this behavior might at first seem pathological, it is in fact quite typical in dynamical systems. Generally speaking, given a group action it is hugely unrealistic to ask for any understanding of every single orbit, since these may typically be arbitrarily complicated. Thus the following result is quite amazing. 

\begin{thm}[Eskin-Mirzakhani-Mohammadi \cite{EM, EMM}]\label{T:main}
The $GL(2,\bR)$ orbit closure of a translation surface is always a manifold. Moreover, the manifolds that occur are locally defined by linear equations in period coordinates. These linear equations have real coefficients and zero constant term.  
\end{thm}

Note that although the local period coordinates are not canonical, if a manifold is cut out by linear equations in one choice of period coordinates, it must also be in any other overlapping choice of period coordinates, because the transition map between these two coordinates is a matrix in $GL(n, \bZ)$. 

Previously orbit closures had been classified in genus 2 by McMullen. (One open problem remains in genus two, which is the classification of $SL(2, \bZ)$ orbits of square-tiled surfaces in $\cH(1,1)$.) The techniques of Eskin-Mirzakhani-Mohammadi's, unlike those of McMullen, are rather abstract, and have surprisingly little to do with translations surfaces. Thus the work of Eskin-Mirzakhani-Mohammadi does not give any information about how many or what sort of submanifolds arise as orbit closures, except for what is given in the theorem statement. 

\section{Applications of Eskin-Mirzakhani-Mohammadi's Theorem}

There are many applications of Theorem \ref{T:main} to translation surfaces, rational billiards, and other related dynamics systems, for example interval exchange transformations. Here we list just a few of the most easily understood applications.

\bold{Generalized diagonals in rational polygons.} Let $P$ be a rational polygon. A generalized diagonal is a billiard trajectory that begins and ends at a corner of $P$. If $P$ is a square, an example is a diagonal of $P$. Let $N_P(L)$ be the number of generalized diagonals in $P$ of length at most $L$. It is a folklore conjecture that 
$$\lim_{L\to\infty}\frac{N_P(L)}{L^2}$$
exists for every $P$ and is non-zero. Previously, Masur had shown that the limsup and liminf are non-zero and finite \cite{Mas, Mas2}. Eskin-Mirzakhani-Mohammadi give the best result to date: with some additional Cesaro type averaging, the conjecture is true, and furthermore only countably many real numbers may occur as such a limit. 

\bold{The illumination problem.} Given a polygon $P$ and two points $x$ and $y$, say that $y$ is illuminated by $x$ if there is a billiard trajectory going from $y$ to $x$. This terminology is motivated by thinking of $P$ as a polygonal room whose walls are mirrors, and thinking of a candle placed at $x$. The light rays travel along billiard trajectories. We emphasize that the polygon need not be convex. 

Leli\`evre, Monteil and Weiss have shown that if $P$ is a rational polygon, for every $x$ there are at most finitely many $y$ not illuminated by $x$ \cite{LMW}. 

\bold{The Wind Tree Model.} This model arose from physics, and is sometimes called the Ehrenfest model. Consider the plane with periodically shaped rectangular barriers (``trees"). Consider a particle (of ``wind") which moves at unit speed and collides elastically with the barriers.  

Delecroix-Hubert-Leli\`evre have determined the divergence rate of the particle for all choices of size of the rectangular barriers \cite{DHL}. Without Eskin-Mirzakhani-Mohammadi (and work of Chaika-Eskin \cite{EC}), the best that could  proven was the existence of a unspecified full measure set of choices of sizes for which such a result holds. 

There are many other examples along these lines, where previously results were known to hold for almost all examples without being known to hold in any particular example, and now with Eskin-Mirzakhani-Mohammadi can be upgraded to hold in all cases. 

\begin{figure}[h!]
\centering
\includegraphics[scale=0.4]{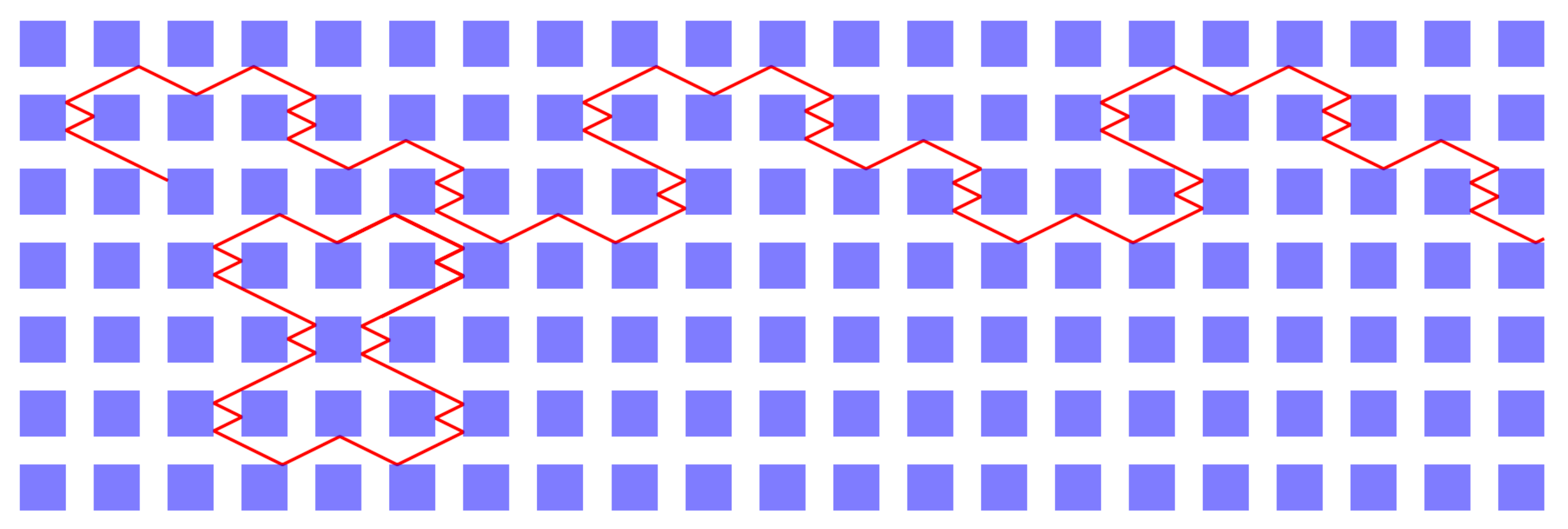}
\caption{An example of a trajectory in the Wind Tree Model. Figure courtesy of Vincent Delecroix.}
\label{F:Vincent}
\end{figure}

\bold{Applications of Eskin-Mirzakhani's proof.} The ideas that Eskin-Mirzakhani developed have applications beyond moduli spaces of translation surfaces. They are currently being used by Rodriguez-Hertz and Brown to study random diffeomorphisms on surfaces \cite{BRH} and the Zimmer program (lattice actions on manifolds), and are also expected to have applications in ergodic theory on homogeneous spaces.

\section{Context from homogeneous spaces}

The primary motivation for Theorem \ref{T:main} is the following theorem. 

\begin{thm}[Ratner's Theorem]
Let $G$ be a Lie group, and let $\Gamma\subset G$ be a  lattice. Let $H\subset G$ be a subgroup generated by unipotent one parameter groups. Then every $H$ orbit closure in $G/\Gamma$ is a manifold, and moreover is a sub-homogeneous space. 
\end{thm}

For example, the theorem applies if $G=SL(3,\bR)$, $\Gamma=SL(3,\bZ)$, and $H=\{h_t: t\in \bR\}$, where
$$h_t=\left(\begin{array}{ccc}1&t&0\\0&1&0\\0&0&1\end{array}\right) \quad\quad\text{or}\quad\quad 
h_t=\left(\begin{array}{ccc}1&t&t^2/2\\0&1&t\\0&0&1\end{array}\right).$$

Ratner's work confirmed a conjecture of Raghunathan. Special cases of this conjecture had previously been verified by Dani and Margulis, for example for the second choice of $h_t$ above \cite{MD}.

The basic idea behind such proofs is the strategy of additional invariance. Given a closed $H$-invariant set, one starts with two points $x$ and $y$ very close together, and applies $h_t$ until the points drift apart. The direction of drift is controlled by another one parameter subgroup, and one tries to show that the closed $H$-invariant set is in fact also invariant under the one parameter group that gives the direction of drift. One continues this argument inductively, each time producing another one parameter group the set is invariant under, until one shows that the closed $H$ invariant set is in fact invariant under a larger group $L$, and is contained in (and hence equal to) a single $L$ orbit. This gives the set in question is homogenous, and in particular a manifold. 

Of course, this is in fact very difficult, and complete proofs of Ratner's Theorem are very long and technical. For one thing, as is the case in the work of Eskin-Mirzakhani, it is in fact too difficult to work directly with closed invariant sets, as we have just suggested. Rather, one first must classify invariant measures. Thus the argument takes place in the realm of ergodic theory, which exactly studies group actions on spaces with invariant measures. See \cite{Morris,Ein} for an introduction. 

The fundamental requirement of the proof is that orbits of nearby points drift apart slowly and in a controlled way. This is intimately tied to the fact that unipotent one parameter groups are polynomial, as can be seen in the $h_t$ above.
Contrast this to the one parameter group $$\left(\begin{array}{cc} e^t&0\\0&e^{-t}\end{array}\right),$$
whose orbit closures may be fractal  sets.   

One might hope to study $GL(2,\bR)$ orbit closures of translation surfaces using the action of 
$$u_t=\left(\begin{array}{cc} 1&t\\0&1\end{array}\right)$$
on strata, in analogy to the proof of Ratner's Theorem. Unfortunately, the dynamics of the $u_t$ action on strata is currently much too poorly understood for this.  %The best known result on the $u_t$ action on strata is the crude, but still very useful, quantitative recurrence result of Minsky-Weiss \cite{MinW}. In particular, this result says that $u_t$ orbits do not diverge.

\section{The structure of the proof} 

The proof of Theorem \ref{T:main}  builds on many ideas in homogeneous space dynamics, including the work of Benoist-Quint \cite{BQ}, and the high and low entropy methods of Einsiedler-Lindenstrauss-Katok \cite{EKL} (see also the beautiful introduction for a general audience by Venkatesh \cite{Venk}). Lindenstrauss won the Fields medal in 2010 partially for the development of the low entropy method, and Benoist and Quint won the Clay prize for their work.

Entropy measures the unpredictability of a system that evolves over time.  

Define $P$ to be the upper triangular subgroup of $SL(2,\bR)$. The proof of Theorem \ref{T:main} proceeds in two main stages. 

In the first,  Eskin-Mirzakhani show that any ergodic $P$  invariant measure is in fact a Lebesgue class measure on a manifold cut out by linear equations, and must be $SL(2,\bR)$ invariant. (An ergodic measure is an invariant measure which is not the average of two other invariant measures in a nontrivial way. Thus the ergodic measures are the building blocks for all other invariant measures.) 

In the second stage, Eskin-Mirzakhani-Mohammadi use this to prove Theorem \ref{T:main}, by constructing a $P$-invariant measure on every $P$-orbit closure. By contrast, it is not possible to directly construct an $SL(2,\bR)$ invariant measure on each $SL(2,\bR)$ orbit closure, and this is why the use of $P$ is crucial. The algebraic structure of $P$ makes it possible to average over larger and larger subsets of $P$ and thus produce $P$ invariant measures, whereas the more complicated algebraic structure of $SL(2,\bR)$ does not allow this.   
(The relevant property is that $P$ is amenable, while $SL(2,\bR)$ is not.) 

In the paper of Eskin-Mirzakhani, which caries out the first stage, the most difficult part is in fact to show $P$-invariant measures are $SL(2,\bR)$ invariant. To do this, extensive entropy arguments are used, partially inspired by the Margulis-Tomanov proof of Ratner's Theorem \cite{MaTo} and to a lesser extent the high and low entropy methods. This part is the technical heart of the argument, and takes almost 100 pages of delicate  arguments. One of the morals is that entropy arguments are surprisingly effective in this context, and can be made to work without the use of an ergodic theorem.  

Once Eskin-Mirzakhani show  that $P$-invariant measures are $SL(2,\bR)$ invariant, they build upon ideas of Benoist-Quint to conclude that every such measure is a Lebesgue class measure on a manifold cut out by linear equations

All together, the proof is remarkably abstract. The only facts used about translation surfaces are formulas for Lyapunov exponents due to Forni \cite{F}. (The Lyapunov exponents of a smooth dynamical system, in this case the action of $g_t$ on a stratum, measure the rates of expansion and contraction in different directions.) Forni was awarded the Brin prize partially for these formulas, which are of a remarkably analytic nature and arose from an insight of Kontsevich \cite{K}.  Eskin-Mirzakhani also make use of a result of Filip \cite{Fi1} to handle a  volume normalization issue at the end of the proof.  

Results which classify invariant measures are rare gems. The arguments are abstract, and their purpose is to rule out nonexistent objects, and thus they cannot be guided by examples. To obtain a truly new measure classification result, one must find truly new ideas from among the sea of ideas which don't quite work, and the devil is in the details. This is very much the case in the paper of Eskin-Mirzakhani.

\section{Relation to Teichm\"uller theory and algebraic geometry}  

Every translation surface has, in particular, the structure of a Riemann surface $X$. The extra flat structure is determined by additionally specifying an Abelian differential $\omega$ (a.k.a. holomorphic one form, or global section of the canonical bundle of $X$). The holomorphic one form is $dz$ on the polygons, where $z$ is the usual coordinate on the plane $\bC\simeq \bR^2$. It has zeros at the singularities of the flat metric. 

Thus every translation surface can be given as a pair $(X,\omega)$. For example, the translation surface given by a square with unit area and opposite edges identified is $(\bC/\bZ[i], dz)$. 

There is a projection map $(X,\omega)\mapsto X$ from a stratum of translation surfaces of genus $g$ to the moduli space $\cM_g$ of Riemann surfaces of genus $g$. Under this map, $g_t$ orbits of translation surfaces project to geodesics for the Teichm\"uller metric. But it is important to note that there is no $GL(2,\bR)$ or $g_t$ action on $\cM_g$ itself, only on (the strata) of the bundle of Abelian differentials over $\cM_g$. This is somewhat analogous to the fact that, given a Riemannian manifold, the geodesic flow is defined on the tangent bundle, and there is no naturally related flow on the manifold itself. 

The map $(X,\omega)\mapsto X$ has fibers of real dimension two (given by multiplying $\omega$ by any complex number), and thus the projection of a four dimensional $GL(2,\bR)$ orbit to $\cM_g$ is a two real dimensional object. It turns out that this object is an isometrically immersed copy of the upper half plane in $\bC$ with its hyperbolic metric: such objects are called complex geodesics or Teichm\"uller disks. Note that Royden showed that the Teichm\"uller metric is equal to the Kobayashi metric on $\cM_g$. 

McMullen showed that every $GL(2,\bR)$ orbit closure in genus 2 is either a closed orbit, or an eigenform locus, or a stratum \cite{Mc}. In particular, every $GL(2,\bR)$ orbit closure of genus 2 translation surfaces is a quasi-projective variety. The corresponding statement for $\cM_2$ is that every complex geodesic is either closed, or dense in a Hilbert modular surface, or dense in $\cM_2$. 

The appearance of algebraic geometry in the study of orbit closures was unexpected, and arose in very different ways from work of McMullen and Kontsevich.  A recent success in this direction is the following, which builds upon Theorem \ref{T:main} and work of M\"oller \cite{M, M2}. 

\begin{thm}[Filip \cite{Fi2}]
In every genus, every $GL(2,\bR)$ orbit closure is an algebraic variety that parameterizes pairs $(X,\omega)$ with special algebro-geometric properties, such as $\Jac(X)$ having real multiplication. 
\end{thm}

Furthermore, Filip's work gives an algebro-geometric characterization of orbit closures: They are loci of $(X,\omega)$ with special algebraic properties, when these loci have maximal possible dimension. (When these loci do not have maximum possible dimension, they do not give orbit closures are not of interest in dynamics.) It is of yet unclear how to make use of this characterization, since calculating the dimension of the relevant sub-varieties of $\cM_g$ seems to be an incredibly difficult problem, very close to the Schottky problem. Nonetheless Filip's work means that an equivalent definition of orbit closure can be given to an algebraic geometer, without mentioning flat geometry or even the $GL(2, \bR)$ action! 

It is at present a major open problem to classify $GL(2,\bR)$ orbit closures. Progress is ongoing using techniques based on flat geometry and dynamics, see for example \cite{Wfield, NW, ANW}. One of the key tools is the author's Cylinder Deformation Theorem, which allows one to produce deformations of a translation surface that remain in the $GL(2, \bR)$ orbit closure, without actually computing any surfaces in the $GL(2, \bR)$ orbit \cite{Wcyl}.

\section{What to read next}

We recommend the two page ``What is \ldots  measure rigidity" article by Einseidler \cite{WhatIs}, as well as the eight page note ``The mathematical work of Maryam Mirzakhani" by McMullen \cite{MonM} and the 13 page note ``The magic wand theorem of A. Eskin and M. Mirzakhani" by Zorich \cite{Z2, Z3}. There are  a large number of surveys on translation surfaces, for example \cite{MT, Z} and the author's recent introduction aimed at a broad audience \cite{Wsurvey}. Alex Eskin has given a mini-course on his paper with Mirzakhani, and notes are available on his website \cite{Enotes}.

\bold{Acknowledgements.} I am very grateful to Yiwei She, Max Engelstein, Weston Ungemach, Aaron Pollack, Amir Mohammadi,  Barak Weiss, Anton Zorich, and Howard Masur for making helpful comments and suggestions, which lead to significant improvements. This expository article grew out of notes for a talk in the Current Events Bulletin at the Joint Meetings in January 2015. I thank David Eisenbud for organizing the Current Events Bulletin and inviting me to participate, and  Susan Friedlander for encouraging me to publish this article. 

%%%%%%%%%%%%%%%%%%%%%%%%%%%%%%%%%%%%%%%%%%
%%%%%%%%%%%%%%%%%%%%%%%%%%%%%%%%%%%%%%%%%%
%%%%%%%%%%%%%%  REFERENCES  %%%%%%%%%%%%%%%%%%
%%%%%%%%%%%%%%%%%%%%%%%%%%%%%%%%%%%%%%%%%%
%%%%%%%%%%%%%%%%%%%%%%%%%%%%%%%%%%%%%%%%%%

\bibliography{mybib}{}

\providecommand{\bysame}{\leavevmode\hbox to3em{\hrulefill}\thinspace}
\providecommand{\MR}{\relax\ifhmode\unskip\space\fi MR }
\providecommand{\MRhref}[2]{%
  \href{http://www.ams.org/mathscinet-getitem?mr=#1}{#2}
}
\providecommand{\href}[2]{#2}
\begin{thebibliography}{McM14}

\bibitem[ANW]{ANW}
David Aulicino, Duc-Manh Nguyen, and Alex Wright, \emph{Classification of
  higher rank orbit closures in $\mathcal{H}^{\rm odd}(4)$}, preprint, arXiv
  1308.5879 (2013), to appear in Geom. Top.

\bibitem[BQ09]{BQ}
Yves Benoist and Jean-Fran{\c{c}}ois Quint, \emph{Mesures stationnaires et
  ferm\'es invariants des espaces homog\`enes}, C. R. Math. Acad. Sci. Paris
  \textbf{347} (2009), no.~1-2, 9--13.

\bibitem[BRH]{BRH}
Aaron Brown and Federico Rodriguez~Hertz, \emph{Measure rigidity for random
  dynamics on surfaces and related skew products}, preprint, arXiv 1406.7201
  (2014).

\bibitem[CE]{EC}
Jon Chaika and Alex Eskin, \emph{Every flat surface is {B}irkhoff and
  {O}seledets generic in almost every direction}, preprint, arXiv 1305.1104
  (2014).

\bibitem[DHL14]{DHL}
Vincent Delecroix, Pascal Hubert, and Samuel Leli{\`e}vre, \emph{Diffusion for
  the periodic wind-tree model}, Ann. Sci. \'Ec. Norm. Sup\'er. (4) \textbf{47}
  (2014), no.~6, 1085--1110.

\bibitem[DM90]{MD}
S.~G. Dani and G.~A. Margulis, \emph{Orbit closures of generic unipotent flows
  on homogeneous spaces of {${\rm SL}(3,{\bf R})$}}, Math. Ann. \textbf{286}
  (1990), no.~1-3, 101--128.

\bibitem[Ein06]{Ein}
Manfred Einsiedler, \emph{Ratner's theorem on {${\rm SL}(2,\Bbb R)$}-invariant
  measures}, Jahresber. Deutsch. Math.-Verein. \textbf{108} (2006), no.~3,
  143--164.

\bibitem[Ein09]{WhatIs}
\bysame, \emph{What is {$\dots$} measure rigidity?}, Notices Amer. Math. Soc.
  \textbf{56} (2009), no.~5, 600--601.

\bibitem[EKL06]{EKL}
Manfred Einsiedler, Anatole Katok, and Elon Lindenstrauss, \emph{Invariant
  measures and the set of exceptions to {L}ittlewood's conjecture}, Ann. of
  Math. (2) \textbf{164} (2006), no.~2, 513--560.

\bibitem[EM]{EM}
Alex Eskin and Maryam Mirzakhani, \emph{Invariant and stationary measures for
  the {$SL(2,\bR)$} action on moduli space}, preprint, arXiv 1302.3320 (2013).

\bibitem[EMM]{EMM}
Alex Eskin, Maryam Mirzakhani, and Amir Mohammadi, \emph{Isolation theorems for
  {$SL(2, \bR)$}-invariant submanifolds in moduli space}, preprint,
  arXiv:1305.3015 (2013).

\bibitem[Esk]{Enotes}
Alex Eskin, \emph{Lectures on the ${SL}(2, \mathbb{R})$ action on moduli
  space}, \url{math.uchicago.edu/~eskin/luminy2012/lectures.pdf}.

\bibitem[Fila]{Fi2}
Simion Filip, \emph{Semisimplicity and rigidity of the {K}ontsevich-{Z}orich
  cocycle}, preprint, arXiv:1307.7314 (2013).

\bibitem[Filb]{Fi1}
\bysame, \emph{Splitting mixed {H}odge structures over affine invariant
  manifolds}, preprint, arXiv:1311.2350 (2013).

\bibitem[For02]{F}
Giovanni Forni, \emph{Deviation of ergodic averages for area-preserving flows
  on surfaces of higher genus}, Ann. of Math. (2) \textbf{155} (2002), no.~1,
  1--103.

\bibitem[KMS86]{KMS}
Steven Kerckhoff, Howard Masur, and John Smillie, \emph{Ergodicity of billiard
  flows and quadratic differentials}, Ann. of Math. (2) \textbf{124} (1986),
  no.~2, 293--311.

\bibitem[Kon97]{K}
M.~Kontsevich, \emph{Lyapunov exponents and {H}odge theory}, The mathematical
  beauty of physics ({S}aclay, 1996), Adv. Ser. Math. Phys., vol.~24, World
  Sci. Publ., River Edge, NJ, 1997, pp.~318--332.

\bibitem[KZ03]{KZ}
Maxim Kontsevich and Anton Zorich, \emph{Connected components of the moduli
  spaces of {A}belian differentials with prescribed singularities}, Invent.
  Math. \textbf{153} (2003), no.~3, 631--678.

\bibitem[LMW]{LMW}
Samuel Leli\`ever, Thierry Monteil, and Barak Weiss, \emph{Everything is
  illuminated}, preprint, arXiv 1407.2975 (2014).

\bibitem[Mas82]{Ma2}
Howard Masur, \emph{Interval exchange transformations and measured foliations},
  Ann. of Math. (2) \textbf{115} (1982), no.~1, 169--200.

\bibitem[Mas88]{Mas2}
\bysame, \emph{Lower bounds for the number of saddle connections and closed
  trajectories of a quadratic differential}, Holomorphic functions and moduli,
  {V}ol.\ {I} ({B}erkeley, {CA}, 1986), Math. Sci. Res. Inst. Publ., vol.~10,
  Springer, New York, 1988, pp.~215--228.

\bibitem[Mas90]{Mas}
\bysame, \emph{The growth rate of trajectories of a quadratic differential},
  Ergodic Theory Dynam. Systems \textbf{10} (1990), no.~1, 151--176.

\bibitem[Mas92]{Ma}
\bysame, \emph{Hausdorff dimension of the set of nonergodic foliations of a
  quadratic differential}, Duke Math. J. \textbf{66} (1992), no.~3, 387--442.

\bibitem[McM03]{Mc}
Curtis~T. McMullen, \emph{Billiards and {T}eichm\"uller curves on {H}ilbert
  modular surfaces}, J. Amer. Math. Soc. \textbf{16} (2003), no.~4, 857--885
  (electronic).

\bibitem[McM14]{MonM}
Curtis McMullen, \emph{The mathematical work of {M}aryam {M}irzakhani},
  Proceedings of the International Congress of Mathematicians (Seoul 2014),
  vol.~1, Kyung Moon SA Co. Ltd., Berlin, 2014, pp.~73--79.

\bibitem[M{\"o}l06a]{M2}
Martin M{\"o}ller, \emph{Periodic points on {V}eech surfaces and the
  {M}ordell-{W}eil group over a {T}eichm\"uller curve}, Invent. Math.
  \textbf{165} (2006), no.~3, 633--649.

\bibitem[M{\"o}l06b]{M}
\bysame, \emph{Variations of {H}odge structures of a {T}eichm\"uller curve}, J.
  Amer. Math. Soc. \textbf{19} (2006), no.~2, 327--344 (electronic).

\bibitem[Mor05]{Morris}
Dave~Witte Morris, \emph{Ratner's theorems on unipotent flows}, Chicago
  Lectures in Mathematics, University of Chicago Press, Chicago, IL, 2005.

\bibitem[MT94]{MaTo}
G.~A. Margulis and G.~M. Tomanov, \emph{Invariant measures for actions of
  unipotent groups over local fields on homogeneous spaces}, Invent. Math.
  \textbf{116} (1994), no.~1-3, 347--392.

\bibitem[MT02]{MT}
Howard Masur and Serge Tabachnikov, \emph{Rational billiards and flat
  structures}, Handbook of dynamical systems, {V}ol.\ 1{A}, North-Holland,
  Amsterdam, 2002, pp.~1015--1089.

\bibitem[NW14]{NW}
Duc-Manh Nguyen and Alex Wright, \emph{Non-{V}eech surfaces in
  {$\mathcal{H}^{\rm hyp}(4)$} are generic}, Geom. Funct. Anal. \textbf{24}
  (2014), no.~4, 1316--1335.

\bibitem[Vee82]{V2}
William~A. Veech, \emph{Gauss measures for transformations on the space of
  interval exchange maps}, Ann. of Math. (2) \textbf{115} (1982), no.~1,
  201--242.

\bibitem[Ven08]{Venk}
Akshay Venkatesh, \emph{The work of {E}insiedler, {K}atok and {L}indenstrauss
  on the {L}ittlewood conjecture}, Bull. Amer. Math. Soc. (N.S.) \textbf{45}
  (2008), no.~1, 117--134.

\bibitem[Wri14]{Wfield}
Alex Wright, \emph{The field of definition of affine invariant submanifolds of
  the moduli space of abelian differentials}, Geom. Topol. \textbf{18} (2014),
  no.~3, 1323--1341.

\bibitem[Wri15a]{Wcyl}
\bysame, \emph{Cylinder deformations in orbit closures of translation
  surfaces}, Geom. Topol. \textbf{19} (2015), no.~1, 413--438.

\bibitem[Wri15b]{Wsurvey}
\bysame, \emph{Translation surfaces and their orbit closures: {A}n introduction
  for a broad audience}, EMS Surv. Math. Sci. \textbf{2} (2015), no.~1,
  63--108.

\bibitem[Zor]{Z3}
Anton Zorich, \emph{The magic wand theorem of {A}. {E}skin and {M}.
  {M}irzakhani}, English translation of paper in Gaz. Math., arXiv 1502.05654
  (2015).

\bibitem[Zor06]{Z}
\bysame, \emph{Flat surfaces}, Frontiers in number theory, physics, and
  geometry. {I}, Springer, Berlin, 2006, pp.~437--583.

\bibitem[Zor14]{Z2}
\bysame, \emph{Le th\'eor\`eme de la baguette magique de {A}. {E}skin et {M}.
  {M}irzakhani}, Gaz. Math. (2014), no.~142, 39--54.

\end{thebibliography}
\bibliographystyle{amsalpha}
%\ann{To keep the capitals in ``WriA" and ``WriB" the bbl file has to be manually updated.}
\end{document}